\documentclass[a4paper,twoside,11pt,reqno]{amsart}
\usepackage{fullpage}
\usepackage[T1]{fontenc}
\usepackage[utf8]{inputenc}

\usepackage{hyperref}
\usepackage[slantedGreeks, partialup, noDcommand]{kpfonts}
\usepackage{graphicx}
\usepackage[mathscr]{euscript}

\hypersetup{ocgcolorlinks=true,allcolors=testc}
\hypersetup{
     colorlinks   = true,
     citecolor    = black
}
\hypersetup{linkcolor=black}
\hypersetup{urlcolor=black}

\newtheorem{theorem}{Theorem}
\newtheorem{lemma}[theorem]{Lemma}
\newtheorem{question}[theorem]{Question}
\newtheorem{proposition}[theorem]{Proposition}
\newtheorem{corollary}[theorem]{Corollary}

\newtheorem{remark}[theorem]{Remark}
\newtheorem{definition}[theorem]{Definition}
\newtheorem*{definition*}{Definition}
\newtheorem{claim}[theorem]{Claim}

\usepackage{dutchcal}

\DeclareMathOperator*{\be}{\mathbb{E}}

\title{Covering the hypercube, the uncertainty principle, \\
and an interpolation formula}

\author{Paata Ivanisvili}
\address{(P.~I.)
Department of Mathematics,
University of California, Irvine\\
Irvine, CA 92697, USA.}
\email{pivanisv@uci.edu} 

\author{Ohad Klein}
\address{(O.~K.) School of Computer Science and Engineering, 
Hebrew University of Jerusalem\\
Jerusalem, Israel}
\email{ohadkel@gmail.com}

\author{Roman Vershynin}
\address{(R.~V.)
Department of Mathematics,
University of California, Irvine\\
Irvine, CA 92697, USA.}
\email{rvershyn@uci.edu}

\begin{document}

\maketitle
\vspace{-3mm}

\begin{abstract}
We show that the minimal number of skewed hyperplanes that cover the hypercube $\{0,1\}^{n}$ is at least $\frac{n}{2}+1$, and there are infinitely many $n$'s when the hypercube can be covered with $n-\log_{2}(n)+1$ skewed hyperplanes.
The minimal covering problems are closely related to the uncertainty principle on the hypercube, where  we also obtain an interpolation formula for multilinear polynomials on $\mathbb{R}^{n}$ of degree less than $\lfloor n/m \rfloor$ by showing that its coefficients corresponding to the largest monomials can be represented as a linear combination of values of the polynomial over the points   $\{0,1\}^{n}$ whose Hamming weights are divisible by $m$. 
\end{abstract}

\bigskip

{\footnotesize

}


\section{Introduction}

\subsection{Covering the hypercube}\label{skew}
How many affine hyperplanes are needed to cover the hypercube $\{-1,1\}^{n}$? Notice that two  affine hyperplanes $x_{1}=-1$ and $x_{1}=1$ cover the hypercube, and clearly this is the minimal number. However, if one requires that the affine hyperplanes are skewed, i.e., $a_{1}x_{1}+\ldots+a_{n}x_{n}+b=0$ with all $a_{1}, \ldots, a_{n}\neq 0$, then the problem becomes challenging\footnote{In what follows we will omit the word affine and we will be referring to such hyperplanes as skewed hyperplanes.}. 

It follows from Littlewood--Offord inequalities that any skewed hyperplane covers at most $n^{-1/2}$ fraction of the points in $\{-1,1\}^{n}$ (up to a universal constant factor), therefore, one needs at least $\Omega(n^{1/2})$ skewed hyperplanes to cover the hypercube.
In \cite{Yehu}, this lower bound was improved to $\Omega(n^{0.51})$, and recently in~\cite{Klein} to $\Omega(n^{2/3} \log(n)^{-4/3})$ by the second named author of the present paper. 

The family of $n+1$ hyperplanes,  $x_{1}+\ldots+x_{n}=2k-n$ for all $k=0, \ldots, n$, covers the hypercube. In fact, if $n$ is even, one can cover with $n$  skewed hyperplanes just by replacing the two hyperplanes corresponding to $k=0$ and $k=n$ in the previous example with one hyperplane $x_{1}+\ldots+x_{n/2}-x_{n/2 -1}-\ldots - x_{n}=0$.
Moreover, it follows from \cite{Alon} that for even $n$, the upper bound $n$ on the minimal cover is also a lower bound if one restricts the covering to the family of ``regular'' hyperplanes, i.e., the ones $\varepsilon_{1}x_{1}+\ldots+\varepsilon_{n}x_{n}+b=0$ with $\varepsilon_{j}=\pm 1$ for all $j=1, \ldots, n$.

Looking at the results for the case of ``regular'' hyperplane cover in \cite{Alon}, one may suspect that in analogy to Littlewood--Offord problem the sharp lower bound on the minimal skew hyperplane cover should be $n$.  Surprisingly, one can cover the hypercube $\{-1,1\}^{5}$ with the following $4$ skewed hyperplanes 
\begin{align*}
&x_{1}+x_{2}+x_{3}+x_{4}+2x_{5}=0,\\
&x_{1}+x_{2}+x_{3}-x_{4}+2x_{5} = 0,\\
&x_{1}+x_{2}+x_{3}+x_{4}-2x_{5}=0,\\
&x_{1}+x_{2}+x_{3}-x_{4}-2x_{5}=0.
\end{align*}
Also the hypercube $\{-1,1\}^{6}$ can be covered with $5$ skewed hyperplanes (see Section~\ref{examp}).  In fact, one can cover the hypercube $\{-1,1\}^{n}$ with $n-\log_{2}(n)+1$ skewed hyperplanes for infinitely many $n$'s. 

\begin{proposition}\label{mth0}
For any integer $m\geq 1$ the hypercube $\{-1,1\}^{2^{m}+m-1}$  can be covered with  $2^{m}$ skewed hyperplanes. 
\end{proposition}

This proposition shows that the skew hyperplane covering problem is genuinely different from  the original ``regular'' problem solved in \cite{Alon}.  
\begin{question}
What is the minimal number of skewed hyperplane cover of the hypercube $\{-1,1\}^{n}$?
\end{question} 
We prove the following lower bound.
\begin{theorem}\label{mth1}
The minimal number of skewed hyperplane cover of $\{-1,1\}^{n}$ is at least $\frac{n}{2}+1$. 
\end{theorem}

There is a close relation between the minimal hyperplane covering problem and the uncertainty principle on the hypercube.
Let $p(x)$ be a polynomial on $\mathbb{R}^{n}$, and let $\mathrm{supp}(p) = \{x \in \mathbb{R}^{n} : p(x) \neq 0\}$.  Under what conditions on $\mathrm{supp}(p) \cap \{-1,1\}^{n}$ and $\mathrm{deg}(p)$ does it follow that $p(x) \equiv 0$ on $\{-1,1\}^{n}$? 

It turns out that the support of a nonzero low degree polynomial cannot be contained in a skewed hyperplane: 
\begin{theorem}[Linial--Radhakrishnan~\cite{Linial}]\label{mth2}
If $\mathrm{deg}(p)<\frac{n}{2}$ and $\mathrm{supp}(p) \cap \{-1,1\}^{n}$ belongs to a skewed hyperplane, then $p(x) \equiv 0$ on $\{-1,1\}^{n}$. 
\end{theorem}


Observe that Theorem~\ref{mth1} follows from Theorem~\ref{mth2}. Indeed, let $H_{1}, \ldots, H_{k+1}$ be a minimal  skew hyperplane cover of $\{-1,1\}^{n}$. If  $H_{j}'s$ are given via equations $\ell_{j}(x) = a_{1j}x_{1}+\ldots+a_{nj}x_{n}+b_{j}=0$, for all $j=1, \ldots, k+1$, then it follows that $p(x) \ell_{k+1}(x) \equiv0$ on $\{-1,1\}^{n}$, where $p(x) = \prod_{j=1}^{k}\ell_{j}(x)$ is a not identically zero polynomial on $\{-1,1\}^{n}$ of degree at most $k$. Hence,  $\mathrm{supp}(p)\cap\{-1,1\}^{n}$ belongs to $H_{k+1}$ and Theorem~\ref{mth2} implies that $k\geq n/2$.

\vskip0.5cm

After the current paper was completed, independently and concurrently the paper \cite{Lisa} appeared on arXiv where Theorem~\ref{mth1} was derived from Theorem~\ref{mth2} proved in  \cite{Linial} (see Lemma 2 in \cite{Linial}). The proof of  Theorem~\ref{mth2} in  \cite{Linial}  in turn is based on either Combinatorial Nullstellensatz or the spectral properties of the Johnson graph (the authors \cite{Linial} attribute the nonsingularity of the Johnson graph to \cite{Gottlieb}). Our proof of Theorem~\ref{mth2}, given in Section~\ref{s: mth2}, is simple and self-contained.

\subsection{An interpolation formula}
In \cite{Alon} sharp lower bound $n$  on the minimal number of ``regular''  hyperplane cover of the $n$ dimensional hypercube (for even $n$) was based on the following technical observation: if a multilinear polynomial $p(x)$ vanishes on all those points of $\{-1,1\}^{n}$ which have even number of $1$'s in its coordinates, and $\mathrm{deg}(p)<n/2$, then $p$ is identically zero (see Lemma 2.1 in \cite{Alon}). This observation suggests that perhaps the coefficients of the multilinear polynomials of small degree can be reconstructed by its values at sparse points of $\{-1,1\}^{n}$. The goal of this section is to obtain such an interpolation formula. 

Recall that any function $f :\{-1,1\}^{n} \mapsto X$, where $X$ is a normed space, has Fourier--Walsh representation 
\begin{equation}\label{eq:fwrep}
f(x) = \sum_{S \subset \{1,\ldots, n\}} \widehat{f}(S) x^{S},
\end{equation}
for some $\widehat{f}(S)\in X$, where $x= (x_{1}, \ldots, x_{n})$, $x^{S} = \prod_{j \in S} x_{j}$  and $x^{\emptyset}=1$. We say that $f$ has degree $\mathrm{deg}(f)$ if $\widehat{f}(S)=0$ for all $S \subset \{1, \ldots, n\}$ with $|S|>\mathrm{deg}(f)$, and there exists a subset $S$ of cardinality $\mathrm{deg}(f)$ such that $\widehat{f}(S)\neq 0$. 

\begin{definition}
For any integer $m>1$ the symbol $W(m)$ denotes the subset of  $\{-1,1\}^{n}$ consisting of all points $x = (x_{1},\ldots, x_{n}) \in \{-1,1\}^{n}$ such that $\#\{ j : x_{j}=-1\}$ is divisible by $m$. 
\end{definition}

In this section we obtain the following interpolation formula
\begin{theorem}\label{interp}
Let $f \colon\{-1,1\}^{n} \mapsto X$ and let $m\geq 2$ be an integer divisible by $2$ such that  
\begin{align}\label{condp}
\mathrm{deg}(f) \leq \frac{n}{m}-\frac{1}{2}.
\end{align}
Then for any $S \subset [n]$ with $|S| = \deg(f)$, there exists a probability measure $d\mu(x)$ supported on $W(m)$ and a sign function $h : W(m) \mapsto \{-1, 1\}$ such that 
\begin{equation}\label{eq:interp}
\widehat{f}(S) = \int_{W(m)} h(x)f(x) d\mu(x)
\end{equation}
Both $d\mu$ and $h$ depend only on $S, m, \mathrm{deg}(f), n$.
\end{theorem}

The next corollary follows from the theorem 
\begin{corollary}\label{interpc}
If $f \colon \{-1,1\}^{n} \mapsto X$ vanishes on a set $W(m)$ for some even integer $m$ satisfying (\ref{condp}), then $f \equiv 0$. 
\end{corollary}

\begin{remark}
When $m=2$, Corollary~\ref{interpc} is the classical result \cite[Lemma 2.1]{Alon}. 
\end{remark}
\begin{remark}
In the proof of Theorem~\ref{interp} both the measure $d\mu$ and $h(x)$ are constructed explicitly. 
\end{remark}

Notice that since $\widehat{f}(S) = \mathbb{E}f(x) x^{S}$ then $\|\widehat{f}(S)\| \leq \max_{x \in \{-1,1\}^{n}} \|f(x)\|$. However, if $f$ has low degree, then $\max_{x \in \{-1,1\}^{n}} \|f(x)\|$ can be replaced by a maximum over sparse family of points of $\{-1,1\}^{n}$ provided that $|S|=\mathrm{deg}(f)$. Indeed,  Theorem~\ref{interp}   gives the following

\begin{corollary}
Let $f:\{-1,1\}^{n} \mapsto X$, where $X$ is a normed space, be a function whose degree satisfies (\ref{condp}), then
$$
\| \widehat{f}(S) \| \leq \max_{x \in W(m)} \|f(x)\|
$$
for all $S \subset \{1, \ldots, n\}$ with $|S| =\mathrm{deg}(f)$. 
\end{corollary}

\section{Proofs}
\subsection{The proof of Theorem~\ref{mth2}}	\label{s: mth2}
Denote $[n]:=\{1,\ldots, n\}$. Every polynomial $p(x)$ of degree $d$, when restricted to $\{-1,1\}^{n}$, can be written as $f(x) = \sum_{|S|\leq d} c_{S} x^{S}$ for some $c_{S} \in \mathbb{R}$. The assumption in Theorem~\ref{mth2} that the support of $f$ is contained within a skewed hyperplane means that 
\begin{equation}\label{eq:linsupp}
    \forall x \in \{-1,1\}^{n} \colon (a_{1}x_{1}+\ldots+a_{n}x_{n}+b) \sum_{|S|\leq d} c_{S} x^{S} = 0,
\end{equation}
where
\begin{equation}\label{eq:skewed}
    \forall i \colon a_i \neq 0.
\end{equation}
When expanding~\eqref{eq:linsupp}, the right hand side means that all coefficients of the monomials $x^{T}$ must vanish, in particular for degree $d+1$ monomials. This means that
\begin{align}\label{sys1}
\sum_{j \in T} a_{j} c_{T\setminus{j}}  =0
\end{align}
for all $T \subseteq [n]$ with $|T|=d+1$. We can view~\eqref{sys1}  as a system of linear equations $Ac=0$
	in $c=(c_S:\;S \subset [n],\, |S| = d)$. It suffices to prove 
	
\begin{lemma}\label{genn}
We have $\mathrm{Ker}(A) =0$ as long as $n\geq 2d+1$. 
\end{lemma}	
The theorem follows from the lemma as follows:
The lemma implies that in~\eqref{eq:linsupp} we have $c_S=0$ for all $S \subset [n]$ with $|S|=d$.
This means that $p(x)$ is in fact a polynomial of degree $d-1$. Since $d$ was merely defined as the degree of $p$ (and assumed to satisfy $2d+1 \leq n$), we can repeat and deduce $p$ is of degree $d-2$ and similarly of degree $0$. Once $n\geq 2$ a hyperplane can not cover the entire hypercube, so $p$ must be the zero polynomial, concluding the proof of Theorem~\ref{mth2}.

\begin{proof}[Proof of Lemma~\ref{genn}]

Note that it is sufficient to prove the lemma only for $n=2d+1$, and it would follow for any $n \ge 2d+1$.
To see this, let $n > 2d+1$ and let $S\subset [n]$ with $|S|=d$; we must show $c_S=0$.
Fix a set $N \subset [n]$ with $S\subset N$ and $|N|=2d+1$, and focus on equations~\eqref{sys1} for $T \subset N$.
From the $n=2d+1$ case we conclude that all involved variables $c_{S'}$ with $S'\subset N$ are $0$, and in particular $c_S=0$.

Next, we prove the lemma by induction on $d$ for the $n=2d+1$ case.

\medskip\noindent\textbf{Base case:} When $d=1$ and $n=3$, Equation~\eqref{sys1} applied on the sets $T=\{2,3\}$, $\{1,3\}$, $\{1,2\}$ yields:
	\begin{align*}
		a_3c_2+a_2c_3 &= 0 \\
		a_1c_3+a_3c_1 &= 0 \\
		a_1c_2+a_2c_1 &= 0,
	\end{align*}
Or in matrix form:
	\[
	\begin{bmatrix}
		0 & a_3 & a_2 \\
		a_3 & 0 & a_1 \\
		a_2 & a_1 & 0 
	\end{bmatrix}
	\begin{bmatrix}
		c_1 \\ c_2 \\ c_3
	\end{bmatrix}
	= 0.
	\]
	The determinant of this $3 \times 3$ matrix equals $2a_1a_2a_3$, which is nonzero by~\eqref{eq:skewed}. We learn that $c_i = 0$ and the lemma follows in this case.
	
	\medskip\noindent\textbf{Inductive step:} Assume that the lemma holds for $d-1$, we prove it for $d$ and $n=2d+1$.
	Let $c$ be a solution to \eqref{sys1}. In order to complete the induction step, we must show that $c=0$.

	Note that the induction hypothesis applied for $d-1$ and $n-2$ (that has $n-2 = 2(d-1)+1$) implies that $\bar{c}=0$ is the unique solution to the system of equations
	\begin{equation}\label{eq:d-1}
		a_{i_1} \bar{c}_{i_2,i_3,\ldots,i_d} + a_{i_2} \bar{c}_{i_1,i_3,\ldots,i_d} + \cdots + a_{i_d} \bar{c}_{i_1,i_2,\ldots,i_{d-1}}=0
		\quad \forall \{i_1,\ldots,i_d\} \subset [n-2],
	\end{equation}
	where $\{i_1, \ldots, i_d\}$ ranges over all subsets of size $d$ of $[n-2]$.

	Fix $j \coloneqq i_{d+1} \in \{n-1,n\}$, and consider all those linear equations~\eqref{sys1} in $c$ that arise for $T=\{i_1,\ldots,i_d,j\}$ where $\{i_1,\ldots,i_d\} \subset [n-2]$:
	\[
		a_{i_1} c_{i_2,i_3,\ldots,i_d,j} + a_{i_2} c_{i_1,i_3,\ldots,i_d,j} +  \cdots + a_{i_d} c_{i_1,i_2,\ldots,i_{d-1},j} 
		+ a_j c_{i_1,i_2,\ldots,i_d}
		= 0	\quad \forall \{i_1,\ldots,i_d\} \subset [n-2].
	\]
	Moving the last term into the right hand side, we get:
	\begin{equation}\label{eq:d-1-j}
	a_{i_1} c_{i_2,i_3,\ldots,i_d,j} + a_{i_2} c_{i_1,i_3,\ldots,i_d,j} +  \cdots + a_{i_d} c_{i_1,i_2,\ldots,i_{d-1},j} 
	= - a_j c_{i_1,i_2,\ldots,i_d}
	\quad \forall \{i_1,\ldots,i_d\} \subset [n-2].
	\end{equation}
	Compare this system with the one in~\eqref{eq:d-1}.
	This is the same system, up to a relabeling of the variables, and different right hand side.
	Let $A_{d-1}$ be the matrix defining the $d-1$ case, then~\eqref{eq:d-1-j} can be written as
	\[
		A_{d-1} \bar{c} = -a_j u,
	\]
	where $u$ is the vector of variables $c_{i_1,i_2,\ldots,i_d}$ and $\bar{c}$ is the vector of variables $c_{i_1,i_2,\ldots,i_d, j}$ for $\{i_1,i_2,\ldots,i_d\} \subset [n-2]$.
	Note that since $\binom{2d-1}{d-1} = \binom{2d-1}{d}$, the matrix $A_{d-1}$ is square, and is hence invertible by the induction hypothesis.
	In particular,
	\[
	\bar{c} = -a_j A_{d-1}^{-1} u.
	\]
	Note that on the right hand side, the only thing that depends on $j$ is $a_j$. That is, both $A_{d-1}$ and $u$ do not depend on whether $j=n-1$ or $j=n$.
	By comparing the two options $j=n-1$, and $j'=n$, we conclude that for all $i_1,i_2,\ldots,i_d \subset [n-2]$,
	\begin{equation}\label{eq:ratio}
		c_{i_1,i_2,\ldots,i_{d-1},j}/ a_j = c_{i_1,i_2,\ldots,i_{d-1},j'}/a_{j'}.
	\end{equation}	
	Note that in~\eqref{eq:ratio}, we can choose the indices $i_1, \ldots, i_d, j', j$ arbitrarily so long as they are distinct.

	\begin{claim}\label{claim:explicit}
		Equation~\eqref{eq:ratio} implies the existence of a single constant $K \in \mathbb{R}$ such that 
		\begin{equation}\label{eq:explicit}
		c_S=K \prod_{i \in S} a_i,
		\end{equation}
		for all $S \subset [n]$ with $|S|=d$.
	\end{claim}
	Plugging the formula for $c_S$ from Claim~\ref{claim:explicit} into the system of equations~\eqref{sys1}, we get that for any $T \subset [n]$ with $|T|=d+1$,
	\[
		0 = \sum_j a_j c_{T \setminus \{j\}} = \sum_{j \in T} a_j \left(K \prod_{i \in T \setminus \{j\}} a_i\right) = \sum_{j \in T} K \prod_{i \in T} a_i = (d+1) K \prod_{i \in T} a_i
	\]
	In~\eqref{eq:skewed}, we assumed $a_i \neq 0$ for all $i$, so it follows that $K=0$. and in particular $c_S = 0$ for all $S \subset [n]$. The induction step and the lemma follows.

	\begin{proof}[Proof of Claim~\ref{claim:explicit}]
        Let $I \subset [n]$ and $J, K \subset [n]\setminus I$ be sets with $|J|=|K|=d-|I|$. We prove by induction on the size of $J$ and $K$ that
		\begin{equation}\label{eq:explicit-induction}
			c_{I \cup J}/\prod_{j \in J} a_{j} = c_{I \cup K}/\prod_{k \in K} a_{k}.
		\end{equation}
		When $J,K$ are of size $d$, and $I=\emptyset$, then we derive~\eqref{eq:explicit} and the proof is complete.

		\noindent\textbf{Base case:} When $J,K$ are of size $0$,~\eqref{eq:explicit-induction} is obvious as $J=K$.

		\noindent\textbf{Inductive step:} Let $J,K\subset [n]$ be of size $\ell\geq 1$ and let $I\subset [n] \setminus (J\cup K)$ be of size $d-\ell$. We prove~\eqref{eq:explicit-induction}. Take $j \in J$ and $k \in K$ and denote $I'=I\cup \{j\}$, $J'=J\setminus \{j\}$ and $K'=K\setminus \{k\}$.
		Then~\eqref{eq:explicit-induction} follows:
		\begin{align*}
		c_{I \cup J} / \prod_{i \in J} a_{i} &\underbrace{=}_{I\cup J = I' \cup J'} \frac{1}{a_j} c_{I' \cup J'} / \prod_{i \in J'} a_{i} \\
		                                     &\underbrace{=}_{\mathrm{induction}} \frac{1}{a_j} c_{I' \cup K'} / \prod_{i \in K'} a_{i} \\
											 &\underbrace{=}_{\eqref{eq:ratio}} \frac{1}{a_k} c_{I \cup K} / \prod_{i \in K'} a_{i} \\
											 &= c_{I \cup K} / \prod_{i \in K} a_{i}.
		\end{align*}
	\end{proof}
\end{proof}

\subsection{The proof of Theorem~\ref{interp}}
Let $d=\deg(f)$ and assume that $S = \{m, 2m, \ldots, dm\}$ (any $S$ with $|S|=d$ is obtained by relabeling of variables).
Recall Equation~\eqref{eq:fwrep} defining the Fourier representation: for any $u = (u_{1}, \ldots, u_{n}) \in \{-1,1\}^{n}$ we have
\begin{align*}
f(u) = \sum_{|S|\leq d} \widehat{f}(S) u^{S} \quad \text{where} \quad u^{S} = \prod_{j \in S} u_{j}.
\end{align*} 
For $(y_{1}, \ldots, y_{d}) \in \{-1,1\}^{d}$, and $(x_{1}, \ldots, x_{n})\in \{-1,1\}^{n}$, we define $y \circ x \in \{-1,1\}^{n}$
by splitting the vector $x$ into disjoint sets of indices $I_{1}, I_{2}, \ldots, I_{d}, I_{\mathrm{extra}}, I_{\mathrm{rest}}$, where 
\begin{align*}
	&I_{1} = (1, \ldots, m),\\
	&I_{2} = (m+1, \ldots, 2m),\\
	&\ldots \\
	&I_{d} = ((d-1)m+1, \ldots, dm),\\
	&I_{\mathrm{extra}} = (md+1,\ldots, md+m/2),\\
	&I_{\mathrm{rest}} = (md+m/2+1, \ldots, n).
\end{align*}
Then, we define 
\[
 y_{j} x_{I_{j}} \mathrel{\mathop:}= (y_{j}x_{(j-1)m+1},y_{j}x_{(j-1)m+2}, \ldots, y_{j}x_{jm}),
\]
and finally,
\begin{align*}
	y \circ x \mathrel{\mathop:}= (&y_{1}x_{I_{1}}, \ldots, y_{d} x_{I_{d}}, x_{I_{\mathrm{extra}}}, x_{I_{\mathrm{rest}}}) \\
	=
	(& y_{1}x_{1},\; y_{1}x_{2}, \ldots,\; y_{1}x_{m}, \\
                  &y_{2}x_{m+1},\; y_{2}x_{m+2}, \ldots,\; y_{2}x_{2m},\\
				  & \qquad\qquad\vdots \\
				  &y_{d} x_{(d-1)m+1},\; y_{d} x_{(d-1)m+2},\ldots,\; y_{d}x_{md}, \\
				  &x_{md+1},\; x_{md+2},\ldots,\; x_{md+m-\frac{m}{2}}, \ldots,\; x_{n})
\end{align*}
Note that the variables in $I_{\mathrm{extra}}$ and $I_{\mathrm{rest}}$ are unchanged.
The coordinates $I_{\mathrm{rest}}$ do not play a role in the proof (and may be empty, e.g. if we have equality in~\eqref{condp}),
but the coordinates $I_{\mathrm{extra}}$ have the important role of ``parity'' in the proof. 

\begin{claim}\label{claim:interp-dist}
There exists a distribution $\mathcal{D}$ for $x$ (on $\{-1,1\}^{n}$) such that 
\begin{equation}\label{eq:interp-dist}
   \widehat{f}(S) = \be_{x \sim \mathcal{D}}\; \be_{y \sim \mathrm{unif}(\{-1,1\}^{d})} \left[ f(y \circ x) \cdot y_{1}\cdots y_{d} \right], 
\end{equation}
and moreover,
\begin{equation}\label{eq:interp-dist-w}
   y \circ x \in W(m)
\end{equation}
for all $y \in \{-1,1\}^{d}$ and all $x \in \mathrm{supp}(\mathcal{D})$. $\mathcal{D}$ depends only on $d,m,n$ but not on $f$.
\end{claim}

Claim~\ref{claim:interp-dist} concludes the proof of~\eqref{eq:interp} by using the sign function $h = y_1 \cdots y_d$ and the measure $d \mu$ depicting the distribution\footnote{Note that here we define $h$ as a function of $y$ while the measure $d\mu$ is of $\{-1,1\}^{n}$. We claim that~\eqref{eq:interp-dist} implies that $y_1 \cdots y_d$ is uniquely determined from $y \circ x$ for all $y \circ x$ having positive probability. To see this, note that formula~\eqref{eq:interp-dist} does not depend on $f$, yet for $\widetilde{f}(z)=z^S$ it has $1$ on the LHS while the RHS is bounded by $1$ by the triangle inequality. This means $\widetilde{f}(y \circ x) = y_1 \cdots y_d$. Consequently,~\eqref{eq:interp} holds with $h=y_1\cdots y_d=\widetilde{f}(y \circ x)$, which is a function of $y \circ x$.}
of $y \circ x$ where $x \sim D$ and $y \sim \mathrm{unif}(\{-1,1\}^{d})$.

\medskip\noindent\textbf{Proof of Claim~\ref{claim:interp-dist}.}
Observe the formula
 \begin{align}\label{iden1}
 \be_{y \sim \mathrm{unif}(\{-1,1\}^{d})} \;  \left[ f(y\circ x) \, y_{1}\ldots y_{d} \right] =\sum_{T \colon \forall j \colon |T\cap I_{j}|=1} \widehat{f}(T) \, x^{T}.
 \end{align}
 To verify~\eqref{iden1}, we expand $f(y\circ x)$ on the left hand side as $\sum_{T} \widehat{f}(T) (y\circ x)^{T}$. The number of times $y_j$ appears in that expression is exactly $|T\cap I_{j}|$. If $T\cap I_{j}$ is empty for some $j$, then $(y\circ x)^{T}$ does not depend on $y_j$, and the multiplication by $y_1\ldots y_d$ in~\eqref{iden1} zeroes out the term $\widehat{f}(T) (y\circ x)^{T}$. Hence relevant terms are only those with $|T \cap I_j| \geq 1$ for all $j=1,\ldots, d$. But since $f$ is of degree $d$ to begin with, we must have $|T \cap I_j| = 1$ for all $j=1,\ldots, d$.
 
\smallskip\noindent\textbf{The distribution $\mathcal{D}$.}
We describe how each chunk $x_{I_j}$ of $x\sim \mathcal{D}$ is drawn. All chunks are drawn independently of the other chunks, except for $x_{I_{\mathrm{extra}}}$ which is chosen last.
\begin{itemize}
	\item $x_{I_j}$ for $j=1, \ldots, d$:
\begin{align*}
	\Pr[x_{I_j} = (z_1, \ldots, z_m)] = \begin{cases}
		1/m & \text{if } z_1 = \cdots = z_m = 1, \\
		\frac{1}{2\binom{m-2}{m/2-1}} & \text{if } z_m=1 \text{ and exactly $m/2$ among } z_1, \ldots, z_{m-1} \text{ are equal } -1, \\
		0 & \text{otherwise}.
	\end{cases}
\end{align*}
Note that the sum over all probabilities is $1$.
\item $x_{I_{\mathrm{rest}}}$:
	we set $x_{I_{\mathrm{rest}}}=(1, 1, \ldots, 1)$ always. 
\item $x_{I_{\mathrm{extra}}}$:
Count the total number $s$ of $-1$'s in all chunks $x_{I_1}, x_{I_2}, \ldots, x_{I_d}$. Define
\begin{align}\label{eq:x_extra}
	x_{I_{\mathrm{extra}}} = \begin{cases}
		(1, \ldots, 1) & \text{if } m | s, \\
		(-1, \ldots, -1) & \text{otherwise}.
	\end{cases}
\end{align}
Observe that necessarily $s$ is divisible by $m/2$, since each choice of $x_{I_j}$ adds either $0$ or $m/2$ to $s$. For this reason,~\eqref{eq:x_extra} defines $x$ is such a way that $x \in W(m)$. Furthermore, for all $y\in \{-1,1\}^{d}$ we have $y\circ x \in W(m)$ essentially because signs do not matter modulo $2$.

Finally, in order to deduce~\eqref{eq:interp-dist} from~\eqref{iden1}, we must check that for all $T\subseteq \{1, \ldots, n\}$ with $\forall j \colon |T\cap I_{j}|=1$ we have 
\begin{align}\label{eq:x_be}
\be_{x\sim \mathcal{D}}[x^T] = \begin{cases}
	1 & \text{if } T = S, \\
	0 & \text{otherwise}.
\end{cases}
\end{align}
The case $T=S$ is immediate, since by design $x_{jm}=1$ for all $j=1, \ldots, d$.

Suppose $T\neq S$. Focus on a particular coordinate $t \in T$ with $m \nmid t$ and let $j\in \{1, \ldots, d\}$ be the index with $t \in I_j$. Since $T \cap I_j = \{t\}$, we have that $x^T$ is $x_t \cdot x^{T\setminus \{t\}}$. When $x\sim \mathcal{D}$, $x^{T\setminus \{t\}}$ is a random variable independent of $x_t$, as we draw the different chunks independently. Hence $\be_{x\sim \mathcal{D}}[x^T] = \be_{x\sim \mathcal{D}}[x_t] \cdot \be_{x\sim \mathcal{D}}[x^{T\setminus \{t\}}]$. In order to deduce~\eqref{eq:x_be} and finish the proof we just need to check that $\be_{x\sim \mathcal{D}}[x_t] = 0$. Indeed, by definition, the probability that $x_t=1$ is $1/m + \frac{\binom{m-2}{m/2}}{2\binom{m-2}{m/2-1}}$, that is, either $x_{I_j}$ is all $1$'s, or we need to choose $m/2$ locations for $-1$'s in $I_j$ out of $I_j \setminus \{t, mj\}$. This probability is $1/2$, validating $\be_{x\sim \mathcal{D}}[x_t] = 0$, concluding the proof.
\end{itemize}

\subsection{The proof of Proposition~\ref{mth0}: how to cover the hypercube efficiently}\label{examp}

Consider $2^{m}$ skewed hyperplanes 
\begin{align*}
\sum_{j=1}^{2^{m}-1}x_{j} + \sum_{j=0}^{m-1}\pm 2^{j}x_{2^{m}+j}=0.
\end{align*}
Since any odd integer $k$, $-(2^{m}-1)\leq k \leq 2^{m}-1$ can be written as a sum $\sum_{j=0}^{m-1}\pm 2^{j}$ for some choice of signs $\pm$, it follows that these hyperplanes cover the cube $\{-1,1\}^{n}$ with $n=2^{m}+m-1$. \qed

\bigskip

There are other examples that are not produced by the construction above. In particular, for $n=6$ the union of the following $5$ skewed hyperplanes 
\begin{align*}
x_1-x_2+2x_3+x_4+x_5+2x_6=0, \\
x_1-x_2+x_3+x_4+x_5-x_6=0, \\ 
x_1-x_2-x_3+2x_4-2x_5+x_6=0, \\
x_1+x_2+x_3+x_4+x_5-x_6=0, \\
x_1-x_2-3x_3+x_4+x_5-x_6=0.
\end{align*}
cover the hypercube $\{-1,1\}^{6}$. 

\section*{Acknowledgement}

We are very grateful to the anonymous referee for the helpful suggestions and insightful comments, which have made the paper clearer and easier to read.
P.I. acknowledges support from NSF grant CAREER-DMS-2152401. O.K. was supported in part by a grant from the Israel Science Foundation (ISF Grant No. 1774/20), and by a grant from the US-Israel Binational Science Foundation and the US National Science Foundation (BSF-NSF Grant No. 2020643).  R.V. acknowledges support from NSF DMS-1954233, NSF DMS-2027299, U.S. Army 76649-CS, and NSF+Simons Research Collaborations on the Mathematical and Scientific Foundations of Deep Learning.

\end{document}